# STRATA OF WEYL GROUPS

G. Lusztig

## Introduction

**0.1.** We fix a root datum with Weyl group $W$. Let $\text{Irr}(W)$ be the set of isomorphism classes of irreducible representations of $W$ over $\mathbf{Q}$.

Let $\mathcal{P}$ be the set of all prime numbers. Let $\bar{\mathcal{P}} = \mathcal{P} \sqcup \{0\}$. For $p \in \mathcal{P}$ let $\mathbf{k}_p$ be an algebraic closure of the field with $p$ elements; let $\mathbf{k}_0 = \mathbf{C}$. For $p \in \bar{\mathcal{P}}$ let $G_p$ be a connected reductive group defined over $\mathbf{k}_p$ associated to the given root datum (then $W$ is the Weyl group of $G_p$).

**0.2.** For $p \in \bar{\mathcal{P}}$ let $\mathcal{U}_p$ be the set of unipotent classes of $G_p$. This is a finite set, see [L76]. The Springer correspondence gives an (injective) map $s_p : \mathcal{U}_p \to \text{Irr}(W)$ (see [S76] for $p = 0$ or $p \gg 0$ and [L81] for any $p$). Let $\text{Irr}_p(W)$ be the image of $s_p$. It is known that $\text{Irr}_0(W) \subset \text{Irr}_p(W)$. Hence there is a unique injective map $\epsilon_p : \mathcal{U}_0 \to \mathcal{U}_p$ such that $s_0(\xi) = s_p(\epsilon_p(\xi)$ for $\xi \in \mathcal{U}_0$. Let $\mathcal{U}_* = \sqcup_{p \in \bar{\mathcal{P}}} \mathcal{U}_p / \sim$ where $\sim$ is the equivalence relation in which $\epsilon_p(\xi) \sim \epsilon_{p'}(\xi)$ for any $\xi \in \mathcal{U}_0$ and any $p, p'$ in $\bar{\mathcal{P}}$. For any $p \in \bar{\mathcal{P}}$ we can regard $\mathcal{U}_p$ as a subset of $\mathcal{U}_*$ and we have

$$\mathcal{U}_* = \cup_{p \in \bar{\mathcal{P}}} \mathcal{U}_p = \cup_{p \in \mathcal{P}} \mathcal{U}_p.$$

Let

$$\text{Irr}_*(W) = \cup_{p \in \bar{\mathcal{P}}} \text{Irr}_p(W) = \cup_{p \in \mathcal{P}} \text{Irr}_p(W) \subset \text{Irr}(W).$$

We have a bijection

(a) $$\mathcal{U}_* \to \text{Irr}_*(W)$$

whose restriction to $\mathcal{U}_p$ is given by $s_p$ for any $p \in \bar{\mathcal{P}}$.

**0.3.** In [L15] we have defined a surjective map $e_p : G_p \to \text{Irr}_*(W)$. The fibres of this map (called the strata of $G_p$) are subsets of $G_p$ which are unions of $G_p$-conjugacy classes of fixed dimension. The set of strata of $G_p$ is denoted by $Str_p$. Note that $Str_p$ is in canonical bijection with $\text{Irr}_*(W)$ and hence also with $\mathcal{U}_*$.

For any group $H$ let $cl(H)$ be the set of conjugacy classes of $H$.

In [L15, 4.1] we have defined a surjective map $\Phi : cl(W) \to \mathcal{U}_*$. Let $\tilde{\Phi} : W \to \mathcal{U}_*$ be the composition $W \to cl(W) \xrightarrow{\Phi} \mathcal{U}_*$ where the first map associates to $w \in W$ the conjugacy class of $w$. The fibres of $\tilde{\Phi}$ are called the strata of $W$. Let $Str(W)$ be the set of strata of $W$; the strata form a partition of $W$ into nonempty unions of conjugacy classes of $W$. Now $Str(W)$ is in obvious bijection with $\mathcal{U}_*$ hence with $Str_p$ for any $p$.

For $w \in W$ let $\mu_w$ be the dimension of the 1-eigenspace of $w$ acting on the reflection representation of $W$. For $C \in cl(W)$ we set $\mu_C = \mu_w$ where $w \in W$.

We recall that an $C \in cl(W)$ is said to be elliptic if $\mu(C) = 0$. Let $cl(W)_{ell}$ be the set of all elliptic $C \in cl(W)$.

From [L11] it is known that the restriction of $\Phi$ to $cl(W)_{ell}$ is injective. Its image is a subset $\mathcal{U}_{bas}$ (*basic unipotent classes*) of $\mathcal{U}_0$ (hence of $\mathcal{U}_*$).

**0.4.** Assume that $G_p$ is almost simple for some/any $p$. Let $C \in cl_{ell}(W)$. In §1 we shall attach to $C$ a connected semisimple group $H_C$ over $\mathbf{k}_2$ (defined up to isomorphism) which we call the modified centralizer of $C$. This group is always a product of symplectic groups over $\mathbf{k}_2$, see Theorem 1.2.

We will show that the stratum $\tilde{\Phi}^{-1}(\Phi(C))$ of $W$ has a remarkable (and so far unexplained) connection to the modified centralizer $H_C$ which makes more precise a property of $\tilde{\Phi}^{-1}(\Phi(C))$ stated in [L14, 0.1(b)]. See Theorem 1.4.

**0.5.** In [L79] a decomposition of $\mathrm{Irr}(W)$ into subsets called families was defined. Let $ce(W)$ be the set of families of $\mathrm{Irr}(W)$. Let $\rho : \mathrm{Irr}_*(W) \to ce(W)$ be the map which to any $E \in \mathrm{Irr}_*(W)$ associates the family of $\mathrm{Irr}(W)$ that contains $E$. (This map is surjective since any family contains a unique special representation and that representation is in $\mathrm{Irr}_*(W)$.)

For $p \in \bar{\mathcal{P}}$ and $c \in ce(W)$ we set $G_{p,c} = \{g \in G_p; e_p(g) \in c\} = e_p^{-1}(c \cap \mathrm{Irr}_*(W))$. The subsets $G_{p,c}$ (for various $c \in ce(W)$) form a partition of $G_p$ into nonempty subsets indexed by $ce(W)$. These subsets are called *superstrata*; they are unions of strata of $G_p$.

In [L15] it was conjectured (and in [C20] it was proved) that any stratum of $G_p$ is locally closed in $G_p$. In [L21] we conjectured that any superstratum of $G_p$ is locally closed in $G_p$. We now further conjecture that

(a) if $p = 0$, any irreducible component of a superstratum of $G_p$ is a rational homology manifold.

In §2 we will give some evidence for (a) by considering the first nontrivial case of it (or rather a Lie algebra version of it).

**0.6.** We consider the surjective map given by the composition

$W \to cl(W) \to \mathcal{U}_* \to \mathrm{Irr}_* W \to ce(W)$.

The first map is the obvious one; the second map is $\Phi$ in 0.3; the third map is 0.2(a); the fourth map is $\rho$ in 0.5. The inverse image of $c \in ce(W)$ under this composition is denoted by $W^c$. The nonempty sets $W^c$ ($c \in ce(W)$) form a partition of $W$ into sets said to be the superstrata of $W$ (they are unions of strata of $W$).

It turns out thar the number of conjugacy classes in a superstratum of $W$ is often (but not always) equal to the cardinal of the corresponding family. For example if $c$ is the family of $W$ of type $E_8$ such that $|c| = 17$ then the corresponding superstratum of $W$ contains exactly 17 conjugacy classes.

**0.7.** In §3 we discuss the question of existence of a semisimple class in a stratum of $G_0$ whose image under $\Phi$ is in $\mathcal{U}_{bas}$.

## 1. The modified centralizer $H_C$

**1.1.** In this section we assume that $G_p$ is simple for some/any $p$ and (except in 1.10) that $C \in cl(W)_{ell}$. We shall attach to $C$ a semisimple group $H_C$ over $\mathbf{k}_2$, which we call the *modified centralizer* of $C$.

For any $p \in \bar{\mathcal{P}}$ let $u_p$ be an element of $\epsilon_p(\Phi(C)) \in \mathcal{U}_p$ (recall that $\Phi(C) \in \mathcal{U}_{bas}$) and let $\tilde{H}_{C,p}$ be the centralizer of $u_p$ in $G_p$ modulo the unipotent radical of its identity component. Let $H_{C,p}$ be the identity component of $\tilde{H}_{C,p}$, a connected reductive group over $\mathbf{k}_p$.

Let $\bar{\mathcal{P}}(C)$ be the set of all $p \in \bar{\mathcal{P}}$ such that $H_{C,p}$ has the same type as $H_{C,0}$. One can verify in the various cases that one of (a),(b),(c) below holds.

(a) $\bar{\mathcal{P}}(C) = \bar{\mathcal{P}}$; we then set $p_C = 0$, $\tilde{H}_C = \tilde{H}_{C,2}, H_C = H_{C,2}$;

(b) $\bar{\mathcal{P}} - \bar{\mathcal{P}}(C)$ consists of a single number (denoted by $p_C$) and $p_C = 2$; we then set $\tilde{H}_C = \tilde{H}_{C,2}, H_C = H_{C,2}$;

(c) $\bar{\mathcal{P}} - \bar{\mathcal{P}}(C)$ consists of a single number (denoted by $p_C$) and $p_C = 3$; we then set $\tilde{H}_C = \tilde{H}_{C,3}, H_C = H_{C,3}$ and we have $H_C = \{1\}$ (viewed as a semisimple group over $\mathbf{k}_2$).

Note that (c) only appears in types $G_2$ and $E_8$ in which case we use the tables in [M80],[LS12].

This completes the definition of $H_C$.

We have the following result.

**Theorem 1.2.** *(a) $H_C$ is product $\prod_{j \in J_C} H_C^j$ where $J_C$ is a finite indexing set and $H_C^j$ ($j \in J_C$) is a simply connected group of type $C_{n_j}$ over $\mathbf{k}_2$ with $n_j \geq 1$ (type $C_1$ means type $A_1$).*

*(b) The conjugation action of $\tilde{H}_C$ on $H_C$ induces the trivial action on $J_C$ with a single exception: one $C$ (denoted $C_{ex}$) in type $E_8$ when $|J_C| = 2$ and $n_j = 1$ for $j \in J_C$.*

**1.3.** For $t \in \mathbf{N}$ let $[0, t] = \{0, 1, 2, \ldots, t\}$; let $E(C) = \prod_{j \in J_C} [0, n_j]$, see 1.2(a). If $C = C_{ex}$ (see 1.2(b)), so that $E(C) = [0, 1] \times [0, 1]$, let $E(C)/\iota$ be the set of orbits of the involution $\iota : E(C) \to E(C)$, $(i_1, i_2) \mapsto (i_2, i_1)$.

The following result describes the sets $\Phi^{-1}(\Phi(C))$.

**Theorem 1.4.** *There exists a map $\zeta : E(C) \to \Phi^{-1}(\Phi(C))$ such that for any $i_* = (i_j)_{j \in J_C} \in E(C)$ we have $\mu_{\zeta(i_*)} = \sum_{j \in J_C} i_j$ and such that (a) or (b) below hold.*

*(a) If $C \neq C_{ex}$ then $\zeta$ is a bijection;*

*(b) if $C = C_{ex}$ then $\zeta$ induces a bijection $E(C)/\iota \to \Phi^{-1}(\Phi(C))$.*
*In particular, in case (a) we have $|\Phi^{-1}(\Phi(C))| = \prod_{j\in J_C}(n_j+1)$ with $n_j$ as in 1.2(a); in case (b) we have $|\Phi^{-1}(\Phi(C))| = 3$.*

**1.5.** In the remainder of this section we prove Theorems 1.2, 1.4. If $G_p$ is of type $A_n$, then 1.2, 1.4, are immediate. There is only one elliptic $C$, the class of a Coxeter element; $\Phi(C)$ is the class of regular unipotent elements in $G_p$ (any $p$); we have $p_C = 0$, $H_C = \{1\}$.

**1.6.** In this subsection we assume that for any $p$, $G_p$ is adjoint of type $C_n, n \geq 2$. Let $\mathcal{X} = \{1,2,3,\dots,2n\}$ and let $\tau : \mathcal{X} \to \mathcal{X}$ be the involution $1 \mapsto 2n, 2 \mapsto 2n-1, \dots, 2n \mapsto 1$. We identify $W$ with the group of permutations of $\mathcal{X}$ that commute with $\tau$. Let $w \in W$. Then $w$ is a product of cycles, say $\nu_i$ cycles of size $i$ for $i = 1, 2, \dots$. Here $\sum_i i\nu_i \leq 2n$ and $\nu_1, \nu_3, \nu_5, \dots$ are even. Moreover for $i$ even we have $\nu_i = \nu'_i + \nu''_i$ where $\nu'_i$ is the number of $i$-cycles which commute with $\tau$; also $\nu''_i$ is even. This defines a bijection between $cl(W)$ and the set $Z'$ consisting of all collections of numbers $\nu_i$ ($i$ odd) and $\nu'_i, \nu''_i$ ($i$ even) as above. Let $Z$ be the set consisting of all collections of numbers $\nu_i$ ($i$ odd) and $\nu_i$ ($i$ even) as above in which to any even $i$ with $\nu_i > 0$ even, we attach a label $\epsilon_i \in \{0,1\}$. It is known that $Z$ parametrizes $\mathcal{U}_* = \mathcal{U}_2$. By [L12] (see also [L15,4.3]), the map $\Phi : cl(W) \to \mathcal{U}_*$ can be identified with the map $Z' \to Z$ (denoted again by $\Phi$) which takes $(\nu_i(i \text{ odd}), \nu'_1, \nu''_i(i \text{ even}))$ to $(\nu_i(i \text{ odd}), \nu_i = \nu'_i + \nu''_i(i \text{ even}), \epsilon_i)$ where for even $i$ with $\nu_i$ even $> 0$, $\epsilon_i$ is 1 if $\nu'_i > 0$ and $\epsilon_i$ is 0 if $\nu'_i = 0$. Now $cl(W)_{ell}$ corresponds to the subset $Z'_{ell}$ of $Z'$ defined by the conditions: $\nu_i = 0$ ($i$ odd), $\nu''_i = 0$ ($i$ even), $\sum_{ieven} i\nu'_i = 2n$. Then $\Phi$ restricts to a bijection $Z'_{ell} \to Z_{bas}$ where $Z_{bas}$ is the subset of $Z$ defined by the conditions $\nu_i = 0$ ($i$ odd), $\sum_i \nu_i = 2n$, $\epsilon_i = 1$ if $i$ is even and $\nu_i > 0$ is even. For $\nu_* = (\nu_i) \in Z_{bas}$, the set $\Phi^{-1}(\nu_*)$ is the subset of $Z'$ defined by the same $\nu_i$ as in $\nu_*$ and such that $\nu'_i > 0$ for any $i$ such that $\nu_i$ is even, $> 0$. An element in this set is determined by assigning to any even $i$ auch that $\nu_i > 0$ an even number $\nu''_i \in [0, \nu_i]$ (if $\nu_i$ is odd) or an even number $\nu''_i \in [0, \nu_i - 1]$ (if $\nu_i$ is even, $> 0$). This gives a bijection

(a) $$\Phi^{-1}(\nu_*) \to \prod_{i\text{even};\nu_i\text{odd}} [0, (\nu_i - 1)/2] \times \prod_{i\text{even};\nu_i>0\text{even}} [0, (\nu_i - 2)/2].$$

We can identify $\nu_*$ with a unipotent class in $G_p$ with $p = 2$ wich when viewed as a class of $Sp_{2n}(\mathbf{k}_2)$ has $\nu_i$ Jordan blocks of size $i$ for any even $i$. The connected centralizer of an element in this unipotent class has reductive part of the form

(b) $$\prod_{ieven;\nu_i\text{odd}} Sp_{\nu_i-1}(\mathbf{k}_2) \times \prod_{ieven;\nu_i>0\text{even}} Sp_{\nu_i-2}(\mathbf{k}_2).$$

This can be extracted from results in [LS12] (which go back to results in [W63]). Since $H_C = H_{C,2}$ for any $C \in cl(W)_{ell}$, we see that Theorems 0.5, 0.7 hold in our case.

For example if $n = 2$, $cl(W)_{ell}$ consists of two elements; one of them has one cycle of size 4 and has $H_C = \{1\}, p_c = 0$; the other has two cycles of size 2 and has $H_C = \{1\}, p_c = 2$.

If $n = 3$, $cl(W)_{ell}$ consists of three elements; one of them has one cycle of size 6 and has $H_C = \{1\}, p_c = 0$; another one has a cycle of size 4 and a cycle of size 2 and has $H_C = \{1\}, p_c = 0$; the third one has three cycles of size 2 and has $H_C$ of type $C_1$, $p_c = 2$.

**1.7.** In this subsection we assume that for any $p$, $G_p$ is adjoint of type $D_n, n \geq 4$. The Weyl group of $G_p$ will be here denoted by $W'$ to distinguish it from the Weyl group $W$ in 1.6; actually $W'$ can be regarded as the subgroup of that $W$ consisting of even permutations. In this case the subset $\mathcal{U}_{bas}$ of $\mathcal{U}_*$ can be identified with the set $Z_{bas}^{ev}$ consisting of all $(\nu_i) \in Z_{bas}$ in 1.2 which satisfy $\sum \nu_i =$ even. The inverse image of $\nu_* \in Z_{bas}^{ev}$ under $\Phi : cl(W') \to \mathcal{U}_*$ again satisfies 1.2(a) (here we use [L12]). We can identify $\nu_*$ with a unipotent class in $G_p$ with $p = 2$ wich when viewed as a class of $SO_{2n}(\mathbf{k}_2)$ has $\nu_i$ Jordan blocks of size $i$ for any even $i$. The connected centralizer of an element in this unipotent class has reductive part as in 1.6(b). This can be again extracted from results in [LS12] (which go back to results in [W63]). Since $H_{C,p} = H_{C,2}$ for any $C \in cl(W')_{ell}$ and any $p$, we see that Theorems 1.2, 1.4 hold in our case.

For example, if $n = 4$, $cl(W')_{ell}$ consists of three elements; one of them has one cycle of size 6 and one cycle of size 2 and has $H_C = \{1\}, p_C = 0$; the other has two cycles of size 4 and has $H_C = \{1\}, p_C = 0$; the third one has four cycles of size 2 and has $H_C$ of type $C_1$, $p_C = 0$.

**1.8.** In this subsection we assume that for any $p$, $G_p$ is adjoint of type $B_n, n \geq 3$. In this case the proof of Theorems 1.2, 1.4 is entirely similar to that in 1.6.

**1.9.** In this subsection we assume that for any $p$, $G_p$ is a simple adjoint group of type $G_2, F_4, E_6, E_7$ or $E_8$. The results in this section are obtained by combining results in [L12] with results in [LS12], [M80]. In each case we give a list with rows of the form $l; C; \gamma; d_n; H$ where:

$C \in cl(W)_{ell}$ (described by the characteristic poynomial of some $w \in C$ in the reflection representation of $W$);

$l$ is the minimum length of an element of $C$;

$\gamma$ is the unipotent class $\Phi(C) \in \mathcal{U}_0$ in the notation of [S85];

$d_n \in \mathrm{Irr}(W)$ has degree $d$ and $n = (l - r)/2$ where $r$ is the dimension of the reflection representation of $W$;

$H$ is the type of $H_C$. (This notation has one ambiguity in type $F_4$.)

Type $G_2$.

$2; \Phi_6; G_2; 1_0; \{1\}$

$4; \Phi_3; G_2(a_1); 2_1; \{1\}$

$6; \Phi_2^2; \tilde{A}_1; 2_2; \{1\}$

In the row starting with 6 we have $p_C = 3$; in the other rows we have $p_C = 0$.

Type $F_4$.
$4; \Phi_{12}; F_4; 1_0; \{1\}$
$6; \Phi_8; F_4(a_1); 4_1; \{1\}$
$8; \Phi_6^2; F_4(a_2); 9_2; \{1\}$
$12; \Phi_4^2; F_4(a_3); 12_4; \{1\}$
$14; \Phi_2^2\Phi_4; C_3(a_1); 16_5; \{1\}$
$16; \Phi_3^2; \tilde{A}_2 + A_1; 6_6; \{1\}$
$10; (\Phi_2^2\Phi_6)'; B_3; 8_3; C_1$
$10; (\Phi_2^2\Phi_6)''; C_3; 8_3; C_1$
$24; \Phi_2^4; A_1 + \tilde{A}_1; 9_{10}; C_1 \times C_1$
In the row starting with $14, 16$ we have $p_C = 2$; in the other rows we have $p_C = 0$.

Type $E_6$.
$6; \Phi_3\Phi_{12}; E_6; 1_0; \{1\}$
$8; \Phi_9; E_6(a_1); 6_1; \{1\}$
$12; \Phi_3\Phi_6^2; A_5 + A_1; 30_3; \{1\}$
$14; \Phi_2^2\Phi_3\Phi_6; A_5; 15_4; C_1$
$24; \Phi_3^3; 2A_2 + A_1; 10_9; C_1$
In each row we have $p_C = 0$.

Type $E_7$.
$7; \Phi_2\Phi_{18}; E_7; 1_0; \{1\}$
$9; \Phi_2\Phi_{14}; E_7(a_1); 7_1; \{1\}$
$11; \Phi_2\Phi_6\Phi_{12}; E_7(a_2); 27_2; \{1\}$
$13; \Phi_2\Phi_6\Phi_{10}; D_6 + A_1; 56_3; \{1\}$
$17; \Phi_2\Phi_4\Phi_8; D_6(a_1) + A_1; 189_5; \{1\}$
$21; \Phi_2\Phi_6^3; D_6(a_2) + A_1; 315_7; \{1\}$
$15; \Phi_2^3\Phi_{10}; D_6; 35_4; C_1$
$23; \Phi_2^3\Phi_6^2; D_6(a_2); 280_8; C_1$
$25; \Phi_2\Phi_3^2\Phi_6; (A_5 + A_1)''; 70_9; \{1\}$
$33; \Phi_2^3\Phi_4^2; A_3 + A_2 + A_1; 210_{13}; C_1$
$31; \Phi_2^5\Phi_6; D_4 + A_1; 84_{12}; C_2$
$63; \Phi_2^7; 4A_1; 15_{28}; C_3$
In each row we have $p_C = 0$.

Type $E_8$.
$8; \Phi_{30}; E_8; 1_0; \{1\}$
$10; \Phi_{24}; E_8(a_1); 8_1; \{1\}$
$12; \Phi_{20}; E_8(a_2); 35_2; \{1\}$
$14; \Phi_6\Phi_{18}; E_7 + A_1; 112_3; \{1\}$
$16; \Phi_{15}; D_8; 210_4; \{1\}$
$18; \Phi_2^2\Phi_{14}; E_7(a_1) + A_1; 560_5; \{1\}$
$20; \Phi_{12}^2; D_8(a_1); 700_6; \{1\}$
$22; \Phi_6^2\Phi_{12}; E_7(a_2) + A_1; 1400_7; \{1\}$

$24; \Phi_{10}^2; A_8; 1400_8; \{1\}$
$26; \Phi_2^2\Phi_6\Phi_{10}; D_7(a_1); 3240_9; \{1\}$
$28; \Phi_3\Phi_9; D_8(a_3); 2240_{10}; \{1\}$
$30; \Phi_8^2; A_7; 1400_{11}; \{1\}$
$34; \Phi_2^2\Phi_4\Phi_8; D_5 + A_2; 4536_{13}; \{1\}$
$40; \Phi_6^4; 2A_4; 4480_{16}; \{1\}$
$16; \Phi_2^2\Phi_{18}; E_7; 84_4; C_1$
$22; \Phi_4^2\Phi_{12}; D_7; 400_7; C_1$
$24; \Phi_2^2\Phi_6\Phi_{12}; E_7(a_2); 1344_8; C_1$
$26; \Phi_3^2\Phi_{12}; E_6 + A_1; 448_9; C_1$
$42; \Phi_2^2\Phi_6^3; A_5 + A_2; 7168_{17}; C_1$
$44; \Phi_3^2\Phi_6^2; A_5 + 2A_1; 3150_{18}; C_1$
$46; \Phi_2^2\Phi_4^2\Phi_6; D_5(a_1) + A_2; 1344_{19}; C_1$
$48; \Phi_5^2; A_4 + A_3; 420_{20}; C_1$
$32; \Phi_2^4\Phi_{10}; D_6; 972_{12}; C_2$
$60; \Phi_4^4; 2A_3; 840_{26}; C_2$
$80; \Phi_3^4; 2A_2 + 2A_1; 175_{36}; C_2$
$44; \Phi_2^4\Phi_6^2; D_6(a_2); 4200_{18}; C_1 \times C_1$
$46; \Phi_2^2\Phi_3^2\Phi_6; (A_5 + A_1)'; 2016_{18}; C_1 \times C_1$
$66; \Phi_2^4\Phi_4^2; A_3 + A_2 + A_1; 1400_{29}; C_1 \times C_1$
$64; \Phi_2^6\Phi_6; D_4 + A_1; 700_{28}; C_3$
$120; \Phi_2^8; 4A_1; 50_{56}; C_4.$

In the rows starting with 26, 34 and ending with $\{1\}$ we have $p_C = 2$; in the row starting with 30 we have $p_C = 3$; in the other rows we have $p_C = 0$.

The row starting with 44 and ending with $C_1 \times C_1$ corresponds to $C = C_{ex}$ in 1.2(b).

We see that Theorem 1.2 holds in our cases. Using also [L15, 4.4-4.8] we see that Theorem 1.4 holds in our cases.

**1.10.** We now consider a stratum of $W$ which does not contain an elliptic conjugacy class. There is a unique conjugacy class $C$ in the stratum such that $\mu_C$ is minimum possible. We can still attach to $C$ a semisimple group $H_C$ over $\mathbf{k}_2$ as follows. We can find a proper parabolic subgroup $\underline{W}$ of $W$ and $\underline{C} \in cl(\underline{W})_{ell}$ such that $\underline{C} \subset C$. For $p \in \bar{\mathcal{P}}$ let $P_p$ be a parabolic subgroup of $G_p$ with reductive quotient $L_p$ whose Weyl group is $\underline{W}$. Let $\bar{L}_p$ be the adjoint group of $L_p$. We can write $\bar{L}_p = \bar{L}_p^1 \times \bar{L}_p^2 \times \ldots \times \bar{L}_p^k$ where $\bar{L}_p^1, \bar{L}_p^2, \ldots, \bar{L}_p^k$ are adjoint simple. We have correspondingly $\underline{W} = \underline{W}^1 \times \underline{W}^2 \ldots \times \underline{W}^k$ (where $\underline{W}^i$ is the Weyl group of $\bar{L}_p^i$) and $\underline{C} = \underline{C}^1 \times \underline{C}^2 \times \ldots \times \underline{C}^k$ where $\underline{C}^i \in cl(\underline{W}^i)_{ell}$. Then $H_{\underline{C}^i}$ are defined and we set $H_C = H_{\underline{C}^1} \times H_{\underline{C}^2} \times \ldots \times H_{\underline{C}^k}$. The analogue of 1.2(a) still holds. The analogue of 1.4 also holds (with $\zeta$ being a bijection).

## 2. Superstrata

**2.1.** In this section we assume that $p = 0$ and that $G = G_0 = Sp_4(\mathbf{C})$. Now

$G$ is a union of 5 strata denoted $G^8, G^6, G^4, G'^4, G^0$. The notation is chosen so that $G^d(d=8,6,0)$ is the union of all conjugacy classes of dimension $d$; $G^4$ is a union of two conjugacy classes of dimension 4 (a unipotent class and a class equal to $(-1)$ times a unipotent class); $G'^4$ is the unique semisimple conjugacy class of dimension 4. Now $G$ is a union of three superstrata, namely $G^8, G^6 \cup G^4 \cup G'^4, G^0$.

**2.2.** Let $\mathfrak{g}$ be the Lie algebra of $G$. In [L15, 6.2] a map $\tilde{e}:\mathfrak{g} \to \mathrm{Irr}(W)$ is defined; its image is a subset of $\mathrm{Irr}_*(W)$. The nonempty fibres of $\tilde{e}$ are said to be the strata of $\mathfrak{g}$. For any $c \in ce(W)$ we set $\mathfrak{g}_c = \{x \in \mathfrak{g}; \tilde{e}(x) \in c\}$. The subsets $\mathfrak{g}_c$ (for various $c \in ce(W)$) form a partition of $\mathfrak{g}$ into nonempty subsets indexed by $ce(W)$. These subsets are said to be the superstrata of $\mathfrak{g}$; they are unions of strata of $\mathfrak{g}$.

Now $\mathfrak{g}$ is a union of 4 strata denoted $\mathfrak{g}^8, \mathfrak{g}^6, \mathfrak{g}^4, \mathfrak{g}^0$. The notation is chosen so that $\mathfrak{g}^d$ is the set of all $\xi \in \mathfrak{g}$ such that the centralizer of $\xi$ has dimension $10-d$. The superstrata of $\mathfrak{g}$ are $\mathfrak{g}^8, \mathfrak{g}^6 \cup \mathfrak{g}^4, \mathfrak{g}^0$. We shall prove the following infinitesimal analogue of 0.6(a) for our $G$.

(a) $\mathfrak{g}^6 \cup \mathfrak{g}^4$ is locally closed in $\mathfrak{g}$ and any irreducible component of it is a rational homology manifold.

(The analogous statements for $\mathfrak{g}^8$ and $\mathfrak{g}^0$ are obvious.)

**2.3.** We can identify $\mathfrak{g}$ with the vector space of all $4 \times 4$ matrices

(a)
$$\begin{pmatrix} A & B & C & D \\ E & F & G & C \\ H & I & -F & -B \\ J & H & -E & -A \end{pmatrix}$$

with entries in $\mathbf{C}$ with the bracket given by the commutator of matrices.

We have $\mathfrak{g}^6 = \mathfrak{X} \cup \mathfrak{Y} \cup \mathfrak{Z}$ where

$\mathfrak{X}$ is the set of all semisimple elements of $\mathfrak{g}$ with eigenvalues $z, z, -z, -z (z \neq 0)$,

$\mathfrak{Y}$ is the set of all semisimple elements of $\mathfrak{g}$ with eigenvalues $z, -z, 0, 0 (z \neq 0)$,

$\mathfrak{Z}$ is the set of all subregular nilpotent elements of $\mathfrak{g}$.

The following result is known. It can be deduced from [IH].

(b) *The closures of $\mathfrak{X}, \mathfrak{Y}$ in $\mathfrak{g}_6$ are smooth, irreducible subvarieties of $\mathfrak{g}$ of dimension 7; they are equal to $\mathfrak{X} \cup \mathfrak{Z}$, $\mathfrak{Y} \cup \mathfrak{Z}$ respectively.*

**2.4.** Note that $\mathfrak{g}^4$ is a nilpotent orbit in $\mathfrak{g}$. Let $N, N'$ in $\mathfrak{g}^4$ be such that $N, N'$ are part of an $sl_2$-triple in $\mathfrak{g}$. Let $\mathfrak{z}$ be the centralizer of $N'$ in $\mathfrak{g}$. We can take

$$N = \begin{pmatrix} 0 & 0 & 0 & 0 \\ 0 & 0 & 0 & 0 \\ 0 & 0 & 0 & 0 \\ 1 & 0 & 0 & 0 \end{pmatrix}, N' = \begin{pmatrix} 0 & 0 & 0 & 1 \\ 0 & 0 & 0 & 0 \\ 0 & 0 & 0 & 0 \\ 0 & 0 & 0 & 0 \end{pmatrix}.$$

The Slodowy slice $Sl = N + \mathfrak{z}$ consists of the matrices

(a)
$$q(b,c,d,f,g,i) = \begin{pmatrix} 0 & b & c & d \\ 0 & f & g & c \\ 0 & i & -f & -b \\ 1 & 0 & 0 & 0 \end{pmatrix}.$$

for various $(b,c,d,f,g,i)\in\mathbf{C}^6$.

We want to describe the set $(\mathfrak{g}^6\cup\mathfrak{g}^4)\cap Sl$.

For $q=q(b,c,d,f,g,i)$ let $\mathcal{Z}_q$ be the centralizer of $q$ in $\mathfrak{g}$. It consists of all elements 2.3(a) satisfying

$$A=0, B=-iE+fH+bJ, C=cJ+gH+fE, D=dJ+cH+bE,$$

$$bH+fI-iF=0, cE+gF+fG=0, cB=bC, bF+cI+dH=fB+iC,$$

$$bG-cF-dE=gB-fC, bE+iG=gI+cH.$$

We can rewrite the last six equations so that they involve $E,F,G,H,I,J$ but not $A,B,C,D$:

(b) $$bH+fI-iF=0,$$

(c) $$cE+gF-fG=0,$$

(d) $$-(ci+bf)E=(bg-cf)H,$$

(e) $$xH+(ci+bf)J=bF+cI,$$

(f) $$-xE+(bg-cf)J=bG-cF,$$

(g) $$bE+iG=gI+cH.$$

(Here we write $x=f^2+ig-d$.)

Let $Y$ (resp. $X$) be the set of all $q=q(b,c,d,f,g,i)$ such that $\dim\mathcal{Z}_q=2$ (resp. $\dim\mathcal{Z}_q\in\{4,6\}$). We have $X\sqcup Y=Sl$.

Let $q=q(b,c,d,f,g,i)$.

If $bg-cf\neq 0$ then from (b),(g) we can express $H,I$ in terms of $F,E,G$ and from (d) we can express $J$ in terms of $F,E,G$. It follows that $\dim\mathcal{Z}_q\leq 3$, so that $q\in Y$.

If $ci+bf\neq 0$ then from (c),(g) we can express $E,G$ in terms $I,F,H$ and from (e) we can express $J$ in terms of $I,F,H$. It follows that $\dim\mathcal{Z}_q\leq 3$, so that $q\in Y$.

In the remainder of the argument we assume that $bg-cf=0$, $ci+bf=0$. Then equation (d) can be omitted and $J$ no longer appears in the equations (e),(f).

If $fx+bc\neq 0$ then from (b),(e) we can express $H,I$ in terms of $F$ and from (c),(f) we can express $E,G$ in terms of $F$. Since $J$ is unrestricted it follows that $\dim\mathcal{Z}_q\leq 2$, so that $q\in Y$.

In the remainder of the argument we assume that $bg-cf=0$, $ci+bf=0$, $fx+bc=0$.

Assume now that $b\ne 0, c\ne 0$. We have $g=cf/b, i=-bf/c$ hence $f^2+ig=0$ and $x=d$, $fd+bc=0$. We have $di-b^2=-bdf/c-b^2=-b(df+bc)/c=0$, $dg+c^2=cdf/b+c^2=c(df+bc)/b=0$, $b^2g+c^2i=bcf-bcf=0$.

It follows that equations (b)-(g) can be rewritten as follows

$H=-(f/b)I+(i/b)F$, $E=(f/c)G-(g/c)F$, $0=0, 0=0, 0=0$.

We see that $F, I, G, J$ are unrestricted and determine $H, E$; thus $\dim \mathcal{Z}_q=4$, so that $q\in X$.

In the remainder of the argument we assume that $bg-cf=0$, $ci+bf=0$, $fx+bc=0, bc=0$.

We have $b(ci+bf)=0$ hence $b^2f=0$ hence $bf=0$. We have $c(bg-cf)=0$ hence $c^2f=0$ hence $cf=0$.

Since $bg-cf=0$, $ci+bf=0$ we must also have $bg=0, ci=0$.

Assume that $c\ne 0$. Then $b=0, f=0, i=0$. The equations (b)-(g) become $cE+gF=0$, $dH=cI$, $-dE=-cF$, $0=gI+cH$. (We omit the two equations $0=0$.) We see that $I=(d/c)H, H=(-g/c)I$, $E=(-g/c)F, F=(d/c)E$ hence $I=(-gd/c^2)I$, $F=(-gd/c^2)F$. If $gd+c^2\ne 0$ it follows that $I=H=E=F=0$ but $G, J$ are unrestricted. Thus $\dim \mathcal{Z}_q=2$ so that $q\in Y$. If $gd+c^2=0$ it follows that $I, F, J, G$ are unrestricted and determine $H, E$; thus $\dim \mathcal{Z}_q=4$ so that $q\in X$.

Assume that $b\ne 0$. Then $c=0, f=0, g=0$. The equations (b)-(g) become $bH-iF=0, dH=bF, -dE=bG, bE+iG=0$. (We omit the two equations $0=0$.) We see that $H=(i/b)F, F=(d/b)H, E=(-i/b)G, G=(-d/b)E$ hence $H=(di/b^2)H$, $E=(di/b^2)E$. If $di-b^2\ne 0$ it follows that $H=E=F=G=0$ but $I, J$ are unrestricted. Thus $\dim \mathcal{Z}_q=2$ so that $q\in Y$. If $di-b^2=0$ it follows that $H, E, J, I$ are unrestricted and determine $F, G$; thus $\dim \mathcal{Z}_q=4$ so that $q\in X$.

We now assume that $b=c=0$. We have $fx=0$. The equations (b)-(g) become $fI-iF=0, gF-fG=0, xH=0, -xE=0, iG=gI$. (We omit the equation $0=0$.)

Assume first that $f\ne 0$ so that $x=0$. The equations above become $I=(i/f)F=0, G=(g/f)F$. (Then $iG=gI$ is automatic.) Now $F, H, E, J$ are unrestricted and determine $I, G$. Thus $\dim \mathcal{Z}_q=4$ so that $q\in X$.

Next we assume that $f=0$. The equations (b)-(g) become $iF=0, gF=0, xH=0, -xE=0, iG=gI$. (We omit the equation $0=0$.)

If $xi\ne 0$, we must have $F=H=E=0, G=(g/i)I$. Now $I, J$ are unrestricted and determine $F, H, E, G$. Thus $\dim \mathcal{Z}_q=2$ and $q\in Y$.

If $xg\ne 0$, we must have $F=H=E=0, I=(i/g)G$. Now $G, J$ are unrestricted and determine $F, H, E, I$. Thus $\dim \mathcal{Z}_q=2$ and $q\in Y$.

We now assume that $xi=0, xg=0$.

Assume first that $x\ne 0$ hence $i=0, g=0$. The equations (b)-(g) become $H=0, E=0$. Now $I, J, F, G$ are unrestricted and $\dim \mathcal{Z}_q=4$ so that $q\in X$.

Assume next that $x=0$ and $i\ne 0$. The equations (b)-(g) become $F=0,G=(g/i)I$. Now $I,J,H,E$ are unrestricted and $\dim \mathcal{Z}_q=4$ so that $q\in X$.

Assume next that $x=0$ and $i=0,g\ne 0$. The equations (b)-(g) become $F=0,I=0$. Now $G,J,H,E$ are unrestricted and $\dim \mathcal{Z}_q=4$ so that $q\in X$.

Assume next that $x=0$ and $i=0,g=0$. The equations (b)-(g) become $0=0$ and $\dim \mathcal{Z}_q=6$ so that $q\in X$.

We see that $X=X_1\sqcup X_2\sqcup X_3\sqcup X_4\sqcup X_5$ where $X_*\subset \mathbf{C}^6$ are described by the following conditions:

$X_1$: $b\ne 0,c\ne 0,bg-cf=0,ci+bf=0,fd+bc=0$;

$X_2$: $b=0,c\ne 0,f=0,i=0,gd+c^2=0$;

$X_3$: $b\ne 0,c=0,f=0,g=0,di-b^2=0$;

$X_4$: $b=0,c=0,f=0,d\ne 0,i=0,g=0$;

$X_5$: $b=0,c=0,f=0,ig-d=0$;

$X_6$: $b=0,c=0,f\ne 0,f^2+ig-d=0$.

Let

$X'=\{q;fd+bc=0,gd+c^2=0,di-b^2=0,d\ne 0\}$,

$X''=\{q;b=0,c=0,f^2+ig-d=0\}$.

We have $X''=X_5\sqcup X_6$, $X'=X'_1\sqcup X'_2\sqcup X'_3\sqcup X'_4$ where

$X'_1=\{q\in X';b\ne 0,c\ne 0\}$

$X'_2=\{q\in X';b=0,c\ne 0\}$

$X'_3=\{q\in X';b=0,c\ne 0\}$

$X'_4=\{p\in X';b=0,c=0\}$.

We show

$X_1=X'_1$.

If $q\in X_1$ we have $g=cf/b,i=-bf/c$ hence $gd+c^2=dcf/b+c^2=c(fd+bc)/b=0$, $di-b^2=-dbf/c-b^2=-b(fd+bc)/c=0$. Also $fd\ne 0$ hence $d\ne 0$ so that $q\in X'_1$. Conversely, if $q\in X'_1$ we have $g=-c^2/d,i=b^2/d$ hence $bg-cf=-bc^2/d-cf=-c(fd+bc)/d=0$, $ci+bf=cb^2/d+bf=b(fd+bc)/d=0$. Thus $q\in X_1$.

We show

$X_2=X'_2$.

If $q\in X_2$ we have $di-b^2=0+0=0$ and $fd+bc=0+0=0$. Also $dg\ne 0$ hence $d\ne 0$ so that $q\in X'_2$. Conversely, if $q\in X'_2$ we have $gd\ne 0,di=0,fd=0$ hence $d\ne 0,i=0,f=0$. Thus $q\in X_2$.

We show

$X_3=X'_3$.

If $q\in X_2$ we have $gd+c^2=0+0=0$ and $fd+bc=0+0=0$. Also $di\ne 0$ hence $d\ne 0$ so that $q\in X'_3$. Conversely, if $q\in X'_3$ we have $di\ne 0,gd=0,fd=0$ hence $d\ne 0,g=0,f=0$. Thus $q\in X_3$.

It is clear that

$X_4=X'_4$.

We see that $X=X'\sqcup X''$.

We now set $\tilde{X}' = \{q; fd+bc=0, gd+c^2=0, di-b^2=0, gi+f^2=0\}$. We have $\tilde{X}'\cap X'' = \{q; b=0, c=0, f^2+ig-d=0, d=0\}$. We have also $X' \subset \tilde{X}'$. Indeed if $q = q(b,c,d,f,g,i) \in X'$ we have $f=-bc/d, g=-c^2/d, i=b^2/d$ hence $f^2+gi=(b^2c^2-b^2c^2)/d^2=0$ as desired. Hence $X'\cup(\tilde{X}'\cap X'')\subset \tilde{X}'$. We show: $\tilde{X}' = X' \sqcup (\tilde{X}'\cap X'')$. Let $q \in \tilde{X}' - (\tilde{X}'\cap X'')$. We must show that $q\in X'$. We have $fd+bc=0, gd+c^2=0, di-b^2=0, gi+f^2=0$. and one of $b=0, c=0, f^2+ig-d=0, d=0$ does not hold. that is one of $b\ne 0, c\ne 0, f^2+ig-d\ne 0, d\ne 0$ holds. If $d\ne 0$ then $q\in X'$. Assume next that $f^2+ig-d\ne 0$. Then $d\ne 0$ and $q\in X'$. Next assume that $b\ne 0$. From $di=b^2$ we have $d\ne 0$ so again $q\in X'$. Next assume that $c\ne 0$. From $gd=c^2$ we have $d\ne 0$ so again $q\in X'$.

We have

$X = \tilde{X}' \cup X''$.

Indeed, $\tilde{X}'\cup X'' = X'\cup(\tilde{X}'\cap X'')\cup X'' = X'\cup X'' = X$. We define $\pi : \mathbf{C}^3 \to \tilde{X}'$ by $(D,G,I)\mapsto q(b,c,d,f,g,i) = q(DI, GD, D^2, -IG, -G^2, I^2)$. We show:

$\pi$ is surjective.

Let $q = q(b,c,d,f,g,i)\in\tilde{X}'$.

If $d=g=i=0$ then $b=c=f=0$ and $\pi^{-1}(q)$ consists of $(0,0,0)$. If $d=g=0, i\ne 0$ we have $b=f=c=0$ and we can find $I\in\mathbf{C}$ such that $I^2=i$. We have $\pi^{-1}(q)=\{(0,0,\epsilon I);\epsilon=\pm 1\}$. If $d=i=0, g\ne 0$ we have $b=f=c=0$ and we can find $G\in\mathbf{C}$ such that $-G^2=g$. We have $\pi^{-1}(q)=\{(0,\epsilon G,0);\epsilon=\pm 1\}$. If $g=i=0, d\ne 0$ we have $b=f=c=0$ and we can find $D\in\mathbf{C}$ such that $D^2=d$. We have $\pi^{-1}(q)=\{(\epsilon D,0,0);\epsilon=\pm 1\}$. If $d=0, i\ne 0, g\ne 0$ we have $b=c=0, f\ne 0$ and we can find $I\in\mathbf{C}$, $G\in\mathbf{C}$ such that $I^2=i, -G^2=g$. We have $I^2G^2=f^2$. Replacing if necessary $G$ by $-G$ we can assume that $-IG=f$. We have $\pi^{-1}(q)=\{(0,\epsilon G,\epsilon I);\epsilon=\pm 1\}$. If $i=0, d\ne 0, g\ne 0$ we have $b=f=0, c\ne 0$ and we can find $D\in\mathbf{C}$, $G\in\mathbf{C}$ such that $D^2=i, -G^2=g$. We have $D^2G^2=c^2$. Replacing if necessary $G$ by $-G$ we can assume that $DG=c$. We have $\pi^{-1}(q)=\{(\epsilon D,\epsilon G,0);\epsilon=\pm 1\}$. If $g=0, d\ne 0, i\ne 0$ we have $c=f=0, b\ne 0$ and we can find $D\in\mathbf{C}$, $I\in\mathbf{C}$ such that $D^2=i, I^2=g$. We have $D^2I^2=b^2$. Replacing if necessary $D$ by $-D$ we can assume that $DI=b$. We have $\pi^{-1}(q)=\{(\epsilon D,0,\epsilon I);\epsilon=\pm 1\}$. If $d\ne 0, i\ne 0, g\ne 0$ we have $b\ne 0, c\ne 0, f\ne 0$ and we can find $D\in\mathbf{C}$, $I\in\mathbf{C}$, $G\in\mathbf{C}$ such that $D^2=d, I^2=i, -G^2=g$. We have $I^2G^2=f^2$. Replacing if necessary $G$ by $-G$ we can assume that $-IG=f$. We have $I^2D^2=b^2$. Replacing if necessary $D$ by $-D$ we can assume that $ID=b$. We have $G^2D^2=c^2$ From $fd+bc=0$ we have $-IGD^2+IDc=0$ so that $c=GD$. We have $\pi^{-1}(q)=\{(\epsilon D,\epsilon G,\epsilon I);\epsilon=\pm 1\}$.

This proves that $\pi$ is surjective and that $\tilde{X}'$ can be identified with the variety of orbits of the $\mathbf{Z}/2$ action on $\mathbf{C}^3$ in which the generator of $\mathbf{Z}/2$ acts by $(D,G,I)\mapsto(-D,-G,-I)$. It follows that $\tilde{X}'$ is irreducible and a rational homology manifold. It is also closed in $\mathbf{C}^6\sim Sl$.

Now $X''$ is an irreducible smooth closed subvariety of $Sl$ isomorphic to $\mathbf{C}^3$. We have $\tilde{X}'\cap X'' = \{q; gi+f^2=0, b=0, c=0, d=0\}$.

We see that $Sl$ has two irreducible components (one smooth, one a rational homology manifold). By properties of Slodowy slices, $\mathfrak{g}^6 \cup \mathfrak{g}^4$ has two irreducible components, one of which is smooth and the other is a rational homology manifold. This completes the proof of 2.2(a).

## 3. Complements

**3.1.** In this section we assume that $p = 0$ and that $G = G_0$ is almost simple.

Let $C \in cl(W)_{ell}$. Let $S_C$ be the stratum of $G$ that contains $\Phi(C) \in \mathcal{U}_0$.

The following result is a consequence of results in [L11a] (at that time strata were not yet introduced).

(a) For $C$ is above (with one exception in type $E_8$) there exists a semisimple class $\gamma \subset S_C$.

Let $(\gamma)$ be the unique piece of $G$ (see [L84,3.1]) that contains $\gamma$ in (a). It is known (see [C20]) that the closure of $(\gamma)$ in $S_C$ is an irreducible component $((\gamma))$ of $S_C$. We have $\dim((\gamma)) = \dim \gamma \dim z_\gamma$ where $z_\gamma$ is the connected centre of the connected centralizer of an element of $\gamma$.

Assume now that $G$ is of type $E_8$. From the tables in [L11,2.3] we see that for some $C \in cl(W)_{ell}$ we can find two semisimple conjugacy classes $\gamma, \gamma'$ both contained in $S_C$ such that $\dim z_\gamma \neq \dim z_{\gamma'}$. Since $\dim \gamma = \dim \gamma'$ it follows that $\dim((g)) \neq \dim((g'))$. We see that $S_C$ may have irreducible components of different dimensions.

*Errata to* [L14]. In the table for $E_8$ on p.64 replace
$(E_7(a_3)J_2)_{E_7A_1}$ by $((D_6A_1)J_2)_{E_7A_1}$ and
$(E_7(a_5)J_2)_{E_7A_1}$ by $((D_6(a_2)A_1)J_2)_{E_7A_1}$.

Also replace $(E_7(a_4)J_2)_{E_7A_1}$ by $((D_6(a_1)A_1)J_2)_{E_7A_1}$ and move it from the row starting with 16 to the row starting with 18.

Department of Mathematics, M.I.T., Cambridge, MA 02139